\documentclass[12pt]{article}

\usepackage{amsmath, amssymb}
\usepackage{graphicx}
\newcommand {\R} {\mathbb{R}}

\newtheorem{Theo}{Theorem}[section]
\newtheorem{Prop}{Proposition}[section]

\newtheorem{Def}{Definition}[section]
\newtheorem{Lema}{Lemma}[section]

\begin{document}
\parindent=1cm
\parskip 1mm

\title {Cubes and adjoints of cross-polytopes}
\author{J. Lawrence \footnote {{\bf Jim Lawrence}, Department of Mathematical Sciences. George Mason University. Fairfax, VA 22030, U.S.A. {\bf   E-mail:} lawrence@gmu.edu}, I.P.  Silva \footnote{ {\bf Ilda P.  Silva}, Departamento de Matemática, Faculdade de Ciencias da Universidade de Lisboa, Campo Grande, 1749-016 Lisboa, PORTUGAL. {\bf E-mail}:  ipsilva@fc.ul.pt}}
\maketitle
\date
{\abstract{We describe a bijection between oriented cubes and adjoints of cross-polytopes. This correspondence is used to prove that the real affine cube is, up to reorientation in the same class,  the unique oriented cube that is realizable. Moreover, its underlying matroid is, up to isomorphism, the unique combinatorial cube that is realizable over  $\R$.}} \vskip 2mm

{\parindent=0cm
{\bf Keywords.} $\pm 1$-vectors, representability over the reals, matroids, oriented matroids, combinatorial cubes, cross-polytopes.\vskip 2mm

{\bf Math. Subsject Classification.} 05B35, 52B40, 52C40, 52C07.}

\section{Introduction}

We consider  $C^n:=\{-1,1\}^n$. The real cube  is the oriented matroid of the affine dependencies of  $C^n$  over  $\R$ and will be denoted ${\bf Q_n}:=Aff(C^n)$.

Oriented cubes  were introduced in \cite {IS1} as a combinatorial approach, via oriented matroids, to linear/affine dependencies  of  $(0,1)-$vectors over the reals. Oriented cubes
 are oriented matroids over  $C^n$  that contain as signed circuits the smallest signed circuits of the real cube - the rectangles - and as signed cocircuits the smallest signed cocircuits of the real cube, the signed cocircuits complementary to the facets and to the skew facets. 

So far the real cube is the unique known example of an oriented cube. The unicity of the real cube would mean a purely combinatorial characterization of the linear and affine dependencies {\it over the reals} of the  $\pm 1$  vectors.

We point out that results of E. Mayhew et al.  \cite{MNW}, see also  \cite {CG}, prove that characterizing real representability/orientability  of matroids in terms of excluded minors would be esssentially useless.

The cross-polytope, $\bf O_n$,  is the oriented matroid of the affine dependencies of  the set   $E=\{e_1,\ldots, e_n\}\cup \{-e_1, \ldots, -e_n\}$  of points of  $\R^n$,  where $e_i:=${\it i-th vector of the canonical basis of} $\R^n$. The cross-polytope is a polytope obtained by a sequence of series extensions of a point with $n$-elements. It is a graphic matroid. In this paper we describe the relationship between combinatorial cubes and adjoints of the cross-polytope.  The main results are: Theorem 3.1  to 3.3. 

Roughly speaking, Theorem 3.1  identifies  oriented cubes and adjoints of the cross-polytope. Theorems 3.2. and 3.3 establish that  the real cube is, up to isomorphism, the unique cube whose underlying matroid is representable over $\R$  and whose orientation is representable over  $\R$.
  
The results are organized in two sections. The next section, section 2, is devoted to general properties of the adjoints and the extension lattice of  the cross-polytope. The main result is Theorem 2.1 which yelds that, up to isomorphism $\bf O_n$ has a unique representable adjoint. A concise description of this canonical adjoint of the cross-polytope is given in Proposition 2.1. Section 3 is used in the proof of the main Theorems. 
 
We assume that the reader is acquainted with matroids and oriented matroids. We suggest  \cite{Ox} as a reference for matroid theory and  \cite{OM}  for oriented matroids.

 We shall be using oriented matroid standard definitions and notation of \cite{OM}. \vskip 5mm
 
 {\parindent=0cm
{\bf Notation.} Given a set  $E$, signed subsets or sign vectors, $X\in \{+,-,0\}^E$, of  $E$   will be denoted as ordered pairs  $X=(X^+,X^-)$ where  $X^+:=\{e\in E: X(e)=+\}$, $X^-:=\{e\in E: X(e)=-\}$  and  $X^0:=\{e\in E: X(e)=0\}$. The support of  $X$  is  $\underline{X}=X^+\cup X^-$  and  $X$  is said to be a signed set complementary to  $X^0$.

For every  $A\subseteq E$,  $X(A)$  denotes the restriction of  $X$ to  $A$.

Given an oriented matroid  ${\cal M}={\cal M}(E)$, its underlying matroid will be denoted  $\underline{\cal M}$.

Given  $E\subset \R^n$ ,   ${\cal L}in(E)$  and  ${\cal A}ff(E)$,  denote, respectively, the oriented matroids of the linear and of the affine dependencies of  $E$  over  $\R$.}

\section{cross-polytopes and their adjoints}

\subsection{The oriented matroid  $\bf O_n$}

The n-cross-polytope, denoted  $\bf O_n$,  is the oriented matroid of affine dependencies of the set  of  $2n$-points of  the real affine space defined by  $E:= \{ e_1, \ldots,\ e_n \}\cup \{-e_1,\ldots,\ -e_n\}$  where  $e_i:=${\it i-th vector of the canonical basis of} $\R^n$.

In order to describe explicitly this oriented matroid  $\bf O_n$  we need some notation.

We denote by  $[n]$  the set   $[n]:=\{1,\ldots, n\}$  and by  $[n]':=\{1',\ldots, n'\}$  a copy of  $[n]$. 
Subsets of  $[n]$  will be denoted  $A,B,\ldots $  and subsets  $[n]'$  by  $A',B',\ldots$.
 
For every subset  $A\subseteq [n]$, the set   $A':=\{i' : i\in A\}$ is the corresponding subset of  $[n]'$.  For every subset  $A'\subseteq [n]'$, the set   $A:=\{i : i'\in A'\}$ is the corresponding subset of  $[n]$.

\begin{Def}
{\rm The n-cross-polytope is the oriented matroid  ${\bf O_n}=O_n(( [n]\cup [n]'), {\cal C})$  whose family of {\it signed circuits} 
${\cal C}$  is defined as:  
$${\cal C}=\{ \pm R_{ij}= ( \{i,j\}, \{i',j'\}) \  1<i<j<n\}.\leqno{(C)}.$$}
\end{Def}

$\bf O_n$  is an oriented matroid of  rank $n+1$. Its
{\it hyperplanes} are the subsets of  $[n]\cup [n]'$  of the following two types:

{\parindent=0cm 
type 1) {\it the LV-facets} $H_A= A\cup ([n]\setminus A)',\ \ A\subseteq [n]$.

type 2) {\it the interior hyperplanes} $H_{ii'}= ([n]\setminus i) \cup ([n]'\setminus i'), \ i=1, \ldots n $.}

The family  $\cal D$  of {\it signed cocircuits} of  ${\bf O_n}$  are the following signed subsets on  $[n]\cup [n]'$, the complements of the two types of hyperplanes:

{\parindent=0cm
type 1)  the positive set  $X_{A'}=(([n]\setminus A)\cup A', \emptyset)$  and its opposite, for every  $A\subset [n]$.

type 2)  $Y_i=(\{i\},\{ i'\})$ and its opposite, for every $i\in [n]$.

$${\cal D}:=\{\pm Y_i\}_{i\in[n]}\cup \{\pm X_{A'}\}_{ A\subseteq [n]}\leqno{(CC)}$$}

Actually, the oriented matroid  ${\bf O_n}$  is a graphic matroid and therefore it has a unique class of orientations \cite{BLV}. \vskip 1cm

\subsection{The canonical adjoint  $\tilde{\bf O}_{\bf n}$  of  $\bf O_n$}

Let  ${\cal M}=(E, {\cal D})$, be an oriented matroid and  $\cal D$  its family of signed cocircuits.  The adjoints and the extension lattice of   $\cal M$  are usually described in terms of  signed vectors  over a  subset  ${\cal D}_1\subseteq {\cal D}$   satisfying the condition:  ${\cal D}={\cal D}_1\uplus -{\cal D}_1$  (see \cite{OM} for the general definitions). 

We describe and consider adjoints and the extension lattice of the cross-polytope {\it with respect to the family  ${\cal D}_1$ } defined by:
$${\cal D}_1:=\{Y_i\}_{i\in[n]}\cup \{X_{A'}\}_{ A\subseteq [n]}$$

An {\it (oriented) adjoint}  of  ${\bf O_n}$,  with respect to  ${\cal D}_1$, is an oriented matroid  $\tilde{\cal M}=\tilde{\cal M}({\cal D}_1)$, over the set  ${\cal D}_1$, satisfying the following two conditions:

{\parindent=0cm

(1) $r(\tilde{\cal M}({\cal D}_1))= r({\bf O_n})=n+1.$

(2) The family  $\tilde{\cal D}$  of signed cocircuits of  $\tilde{\cal M}$  contains the following  $2n$  signed subsets of  ${\cal D}_1$:

$\tilde{\cal X}[i]=(\{ X_{A'}\}_{A\subseteq [n]\setminus i}\cup Y_i,\emptyset),\ \forall\ i\in[n]$  and

$\tilde{\cal X}[i']=(\{ X_{A'\cup i'}\}_{A\subseteq [n]\setminus i}, Y_i),\ \forall\ i'\in[n]'$.} 
\vskip 5mm

{\parindent=0cm
{\bf Remark.} In what follows we shall denote  by $\tilde{Y}$, $\tilde{H}_i$, $\tilde{H}_{i'}$  the following subsets of  ${\cal D}_1$:

$\tilde{Y}:=\{Y_1,\ldots, Y_n\}$,  $\tilde{H}_i:={\cal X}[i]^0= \{ X_{A'+i'}\}_{ A\subseteq [n]\setminus i}\cup (\tilde{Y}\setminus Y_i) $  and  
$\tilde{H}_{i'}:={\cal X}[i']^0= \{ X_{A'}\}_{ A\subseteq [n]\setminus i}\cup (\tilde{Y}\setminus Y_i) $.

By definition of adjoint   $\tilde{H}_i$  and  $\tilde{H}_{i'}$  are hyperplanes of every adjoint.  Actually (Proposition 2.2.2) $\tilde{Y}$  is also a hyperplane of every adjoint of   $\bf O_n$.}
\vskip 5mm

The compositions  of signed cocircuits  $\pm \tilde{X}[e], \ e\in [n]\cup [n']$  are the {\it principal covectors of an(any) adjoint of  $\bf O_n$}. 

There are oriented matroids with no adjoint, however every realizable oriented matroid has at least one realizable adjoint.

Notice that in general, different realizations may produce non-isomorphic realizable adjoints.

We use the canonical realization of the cross-polytope to define  a canonical adjoint of the cross-polytope. We then prove, Theorem 2.1, that, up to isomorphism, this is the unique realizable adjoint of  $\bf O_n$.  
\vskip 5mm

\begin{Def}{( Canonical representation of  $\bf O_n$  and canonical adjoint  $\tilde{\bf O}_{\bf n}$)}

{\rm 
Let  $E_n=(e_1,\ldots, e_n)$  be the canonical basis of  $\R^n$  and consider the subsets  $V_n:=\{ (e_i, 1)\}_{i\in [n]}\ \ and  \ \ V'_n:=\{ (-e_i, 1)\}_{i\in [n]}$  of  $\R^{n+1}$.   
Clearly 
$${\bf O_n}:= Aff (E_n\cup -E_n)\simeq Lin(V_n\cup V'_n)$$

The two types of hyperplanes of the oriented matroid  $Lin(V\cup V')$  are defined analytically by:

{\parindent=0cm
1) For every  $i=1,\ldots, n$ ,  $H_{ii'}:  (e_i,0)\cdot  (x,x_{n+1})=0$.

2) For every  $A\subseteq [n]$, denoting by  $v_A$  the vector of  $\R^n$  defined by  $v_A(i)=-1 \ for \ i\in A$  and  $v_A(i)=1\ for\ i\in [n]\setminus A$,

$H_A:  (v_A,1)\cdot (x,x_{n+1})=0$.} 

Consider  $\tilde{V}:=\{(e_i,0)\}_{i\in [n]} \cup \{(v_A,1)\}_{A\subseteq [n]}\subseteq \R^{n+1}$.  Let   $f: \tilde{V} \longrightarrow  {\cal D}_1$  be the bijection defined by  $f( (e_i,0))=Y_i,\ \forall i\in [n]$  and  $f( (v_A,0))=  X_{A'},\ \forall A'\subseteq[n]$.  The {\it canonical adjoint  $\tilde {\bf O_n}= \tilde{\bf O}_{\bf n}({\cal D}_1)$}, of  $\bf O_n$  with respect to  ${\cal D}_1$,  is the oriented matroid over  ${\cal D}_1$  obtained as  image of  $Lin(\tilde{V})$  by the correspondence  $f$. }
\end{Def}
\vskip 5mm

The next Proposition gives a useful characterizations of the canonical adjoint that will be used latter.\vskip 2mm

\begin{Prop} {\rm (Affine representations of   $\tilde{\bf O}_{\bf n}$  of  ${\bf O_n}$)}  Let  $\tilde{\bf O}$  be the canonical adjoint of the cross-polytope  $\bf O_n$  with respect to  ${\cal D}_1$.

Let  $B=(y_1,\ldots, y_n, b)$  be an affine basis of  $\R^n$  with a distinguished vertex  $b$. For every subset  $A\subseteq [n]$  denote by  $b_A:=\frac {1}{|A|+1}(b+\sum_{i\in A}y_i)$   the barycenter of the face  $F_A:=conv(b\cup \{ y_i\}_{i\in A})$  of the simplex  $conv(B)$.

{\parindent=0.5mm
1. The canonical adjoint  $\tilde{\bf O_n}$  of  $\bf O_n$  is isomorphic, via the correspondence   $f:{\cal D}_1\longrightarrow \R^n$  defined by  $f(Y_i)=y_i$,  $f(X_{A'})=b_{[n]\setminus A}$,  to the acyclic oriented matroid  $Aff(\{y_i\}_{i\in[n]}\cup \{b_A\}_{A\subseteq [n]})$.

2.  For every  $c\in int(conv(B))$   there is a unique affine representation  $f_c:{\cal D}_1\longrightarrow \R^n$  of  $\tilde{\bf O_n}$  satisfying the following condition:   
$$f_c(y_i)=y_i,\ \ f_c(X_{[n]'})=b\  \  and\  \ f_c(X_{\emptyset})=c.$$
The representation  $f_c$  is extended to  every  $X_{A'}$, $A\subseteq [n]$,   defining  $f_c(X_{A´})=S_A\cap  F_{[n]\setminus A}$, where  $F_{[n]\setminus A}:=conv (\{y_i\}_{i\in [n]\setminus A}\cup b)$,   and   $S_A$  is the affine space of  $\R ^n$  defined  by  $S_A:=aff(\{y_i\}_{i\in A}\cup c)$.}

We will denote by  $B(Y_1,\ldots, Y_n, X_{[n]})$  any affine representation of   $\tilde{\bf O_n}$  defined by an affine basis  $B=(y_1,\ldots y_n,b)$  and any choice of  the interior point  $c\in int(conv(B))$  that represents  $X_\emptyset$
\end{Prop}
  
{\parindent=0cm
{\bf Proof.}  Consider the canonical representation of  $\tilde{\bf O_n}\simeq Lin(\tilde{V})$  where   $\tilde{V}:=\{(e_i,0)\}_{i\in [n]} \cup \{(v_A,1)\}_{A\subseteq [n]}\subseteq \R^{n+1}$.

1) To obtain a represention of  $\tilde{\bf O}_{\bf n}$  described in 1), consider the hyperplane  $\tilde{H}$  of  $\R^{n+1}$  defined by:
$${H}: x_1+\ldots+x_n+(n+2)x_{n+1}=2.$$   
For every  $v\in \tilde{ V}$  let  $v':= \lambda^+ v \cap  \tilde{H}$  be the intersection of the positive halfline  $\lambda^+ v ,\ \lambda^+\in \R^+$  with the hyperplane   ${H}$. Consider   
$$\tilde{V}':=\{v'\in {H}:\  v\in{\cal V}\}.$$
By definition of  $\tilde{V}'$,   $\tilde{\bf O}_{\bf n}:= Lin(\tilde{V}) \simeq Aff(\tilde{ V}')$.

We leave it to the reader to verify that the points  $y_i=(2e_i,0)= \lambda^+(e_i,0)\cap {H}$  and  $b=v'_{[n]}:=(v_{[n]},1)=(-1_{[n]},1)$  form an affine basis of the hyperplane $(H)$  and that  for every  $A\subset [n]$,  the halfline  $\lambda^+(v_A,1)$  intersects the hyperplane  ${H}$  in the point  $b_A =\frac{1}{n-|A|+1} (b+\sum_{i\in [n]\setminus A}y_i)$, yelding the result.\vskip 2mm

2) To obtain a representation of   $\tilde{\bf O_n}$  with  $c=\beta b +{\sum _{i=1}}^n \alpha_i y_i$  with    $ \beta, \alpha_i>0$  and   $\beta + {\sum _{i=1}}^n \alpha_i=1$,   consider the hyperplane    
$$H_c: \  \alpha_1 x_1+\ldots+\alpha_n x_n+(\beta+1)x_{n+1}=2\beta $$
As before it is straightforward to prove  that for every  $v\in {\cal V}$  the positive halfline  $\lambda^+ v$  intersects   $H_c$.  Taking  $B=(y_1,\ldots, y_n, b)$  with  $y_i=\frac{2\beta}{\alpha_i} e_i, \ i=1,\ldots, n$ , $b=(v_{[n]},0)$   and   $c=\beta b +{\sum _{i=1}}^n \alpha_i y_i= \beta (v_\emptyset , 1)$. We conclude from  1) that the point  $b_{c_A}:=\lambda^+(v_A,1) \cap H_c$  that satisfies  $b_{c_A}:=\lambda^+(v_A,1) \cap H_c=S_A\cap   F_{[n]\setminus A}$.
}
 \vskip 5mm

\subsection{Extension lattice and properties of the adjoints of ${\bf O_n}$}\vskip 5mm

A {\it single element extension of a matroid}  ${\cal M}={\cal M}(E)$  is a matroid  ${\cal M}\cup p$, over the set  $E\cup p$, with the same rank as $\cal M$  such that  $({\cal M}\cup p)\setminus p={\cal M}$.

Single element extensions of an oriented matroid are described by  {\it localizations}  in  ${\cal M}$.  If  $\cal D$  is the family of signed cocircuits of  $\cal M$  and  ${\cal D}_1\subset \cal D$  is such that  ${\cal D}={\cal D}_1 \uplus -{\cal D}_1$,  a localization in  $\cal M$  with respect to  ${\cal D}_1$  is a  signed subset  $\tilde{\cal E}= (\tilde{\cal E}^+, \tilde{\cal E}^-)\in \{+,-,0\}^{{\cal D}_1}$  that defines the position of the new point  $p$  with respect to the hyperplanes of  $\cal M$ in the following way (for further details see \cite{OM}) :  

{\parindent=0cm

For every  $X=(X^+,X^-)\in \tilde{\cal E}^+$,  $(X^+\cup p, X^-)$  is a signed cocircuit of  ${\cal M}\cup p$. 

For every  $X=(X^+,X^-)\in \tilde{\cal E}^-$,  $(X^+, X^-\cup p)$  is a signed cocircuit of  ${\cal M}\cup p$. 

For every  $X=(X^+,X^-)\in \tilde{\cal E}^0$,  $(X^+, X^-)$  is a signed cocircuit of  ${\cal M}\cup p$.}

A sign vector  $\tilde{\cal E}= (\tilde{\cal E}^+, \tilde{\cal E}^-)\in \{+,-,0\}^{{\cal D}_1}$  defines a localization in  $\cal M$  if it defines an extension in every contraction  ${\cal M}/L$  by a hyperline  $L$.  

We recall that given an oriented matroid  ${\cal M}=(E, {\cal D})$,  $\cal D$  being its family of signed cocircuits, the contraction of  $\cal M$  by  $L$  is the (representable) oriented matroid of rank 2,   ${\cal M}/L =(E\setminus L, {\cal D}/L)$,  with family of signed cocircuits   ${\cal D}/L:=\{X=(X^+,X^-)\in {\cal D} :\ L\subset X^0\}$. A sign vector  $\tilde{\cal E}= (\tilde{\cal E}^+, \tilde{\cal E}^-)\in \{+,-,0\}^{{\cal D}_1}$  defines a localization in  $\cal M$ if and only if  it satisfies the following condition:

{\parindent=0cm
(Loc) {\it For every hyperline  $L$  of  $\cal M$  the restriction  $\tilde{\cal E}/L=(\tilde{\cal E}/L^+, \tilde{\cal E}/L^-)$  of  $\cal E$  to ${\cal D}/L$  is orthogonal to every signed triplet of cocircuits of  ${\cal D}/L$  of the following forms:    $(\{X,Y\}, Z)$  such that  $X^+\cap Y^+\subseteq Z^+$    and   $(\{X,Y,Z\}, \emptyset)$   such that  $X^+\cap Y^+\subseteq Z^-$.}
}

The {\it extension lattice }   $\tilde{\cal E}({\cal M})$  of a matroid  $\cal M$  is the poset of localizations with respect to  ${\cal D}_1$  ordered by the usual order  $\tilde{\cal E}_1= (\tilde{\cal E}_1^+, \tilde{\cal E}_1^-) \preceq \tilde{\cal E}_2= (\tilde{\cal E}_2^+, \tilde{\cal E}_2^-)$  iff  $\tilde{\cal E}_1^+\subseteq \tilde{\cal E}_2^+$  and  $\tilde{\cal E}_1^-\subseteq \tilde{\cal E}_2^-$, endowed with an element  $\hat{1}$.

The next Theorem A of  (Theorem 7.5.8  \cite{OM}  page 311) will be used in the sequel to obtain properties of the adjoints of the cross-polytope.
\vskip 5mm

{\parindent=0cm
{\bf Theorem A} ( Bachem and Kern \cite{BK}) 

{\it Consider an oriented matroid  ${\cal M}=(E,{\cal D})$,  $\cal D$  its family of signed cocicuits.  Let  ${\cal D}_1\subset {\cal D}$  be such that  ${\cal D}={\cal D}_1\cup -{\cal D}_1$.

If   $\tilde{\cal M}=({\cal D}_1, \tilde{\cal L})$  is an adjoint of  $\cal M$,  with  family of signed covectors   $ \tilde{\cal L}$ then every  signed covector  $\tilde{V}= (\tilde{V}^+, \tilde{V}^-)\in \tilde{\cal L}$  is a localization in  $\mathcal M$.  }
}
\vskip 5mm

\begin{Prop}{\rm (Properties of every adjoint of  ${\bf O_n}$)}

Consider the cross-polytope  ${\bf O_n}$  and let  ${\cal D}_1$  be the subfamily of its signed cocircuits  defined by 
$${\cal D}_1:=\{Y_i\}_{i\in[n]}\cup \{X_{A'}\}_{ A\subseteq [n]}.$$
Let  $\tilde{\cal M}=\tilde{\cal M}({\cal D}_1)$  be an adjoint of  ${\bf O_n}$  with respect to  ${\cal D}_1$. Then,
\begin{enumerate}
\item For every  $i\in [n]$ and  every  $A\subseteq [n]\setminus i$ the signed set    
$\tilde{C}(i,A)= (\{Y_i,X_{A'\cup i'}\}, X_{A'})$  is a signed circuit of  $\tilde{\cal M}$.

\item $\tilde{\cal M}$  is an acyclic oriented matroid whose LV-face lattice is isomorphic to the face lattice of a simplex with vertices  $Y_1, \ldots, Y_n, X_{[n]'}$ and whose facets are the hyperplanes  $\tilde{Y}= \{Y_1,\ldots, Y_n\}$  and  $\tilde{\cal H}_{i}:=X[i]^0$  for every  $i\in [n]$.

\item Let  $\tilde{X}[0]$  denote the set of positive signed cocircuits of  $\tilde{\cal M}$  whose complement is the facet  $\tilde{Y}$. $\tilde{X}[0]:=( \{X_{A'}\}_{ A\subseteq [n]}, \emptyset)$  defines a single element extension  ${\bf O_n}\cup 0$  of  ${\bf O_n}$. The oriented matroids  ${\bf O_n}\cup 0$  and  ${\bf O_n}$  have exactly the same number of hyperplanes and of signed cocircuits.

\item For every  $i\in [n]$,  ${\bf O_n}/i \simeq {\bf O_n}/i' \simeq  {\bf O_{n-1}}\cup 0$  where   ${\bf O_{n-1}}\cup 0$  is the extension of  $\bf O_{n-1}$  defined in 3). 

Moreover,

$\tilde{\cal M}(\tilde{H}_i)$  is an adjoint of   ${\bf O_n}/i$  with respect to  ${\cal D}_1/i:=\{X=(X^+,X^-)\in {\cal D}_1:\ i\in X^0\}$  and 

$\tilde{\cal M}(\tilde{H}_{i'})$  is an adjoint of   ${\bf O_n}/i'$  with respect to  ${\cal D}_1/i':=\{X=(X^+,X^-)\in {\cal D}_1:\ i'\in X^0\}$ .
 \end{enumerate}
\end{Prop}

{\parindent=0cm

{\bf Proof.} Let  $\tilde{\cal M}=\tilde{\cal M}({\cal D}_1)$  be an adjoint of the cross-polytope  ${\bf O_n}$.

{\bf 1) }The existence of the rank preserving embedding of  the poset,  $L^{op}$, the poset of flats of  ${\bf O_n}$  with the reverse order, into the lattice of flats of  $\tilde{\cal M}$, implies that for every  $i\in [n]$ and  every  $A\subseteq [n]\setminus i$ the set   $C=\{Y_i, X_{A'+i}, X_{A'}\}$  is a flat and a circuit of rank 2  of  $\tilde{\cal M}$. 
}

On the other hand, by Theorem A every covector  of  $\tilde{\cal M}$  is a localization  in  ${\bf O_n}$  and, by definition of localization with respect to  ${\cal D}_1$,  must be orthogonal to the signature:  $\tilde{C}(i,A)=(\{Y_i, X_{A'+i}\}, X_{A'})$  of  $C$, implying that  $\tilde{C}(i,A)$  must be one of the two opposite signed circuits of  $\cal M$  with support  $C$.

For every  $i\in [n]$ and  every  $A\subseteq [n]\setminus i$ the  three signed cocircuits  $Y_i,X_{A'},\ X_{A'+i'} $  are the signed cocircuits of the contraction  ${\bf O_n}/L$  of  $\bf O_n$  by the hyperline  $L=H_{ii'}\cap H_A \cap H_{A+i}= A\cup ([n]\setminus (A+i))'$, proving 1). \vskip 2mm

{\parindent=0cm
{\bf 2)}  By definition of adjoint for every  $i\in [n]$,  $\tilde{X}[i]:=( \{X_{A'}\}_{A\subseteq[n]\setminus i} \cup (\{Y_i\}, \emptyset)$  is a positive signed cocircuit of  $\tilde{\cal M}$  and consequently  $\tilde{H}_i:=\tilde{X}[i]^0$ is a LV-facet of  $\tilde{\cal M}$.
}

The circuits   $\tilde{C}(i,A)$  of  $\tilde{\cal M}$  imply that if  $\tilde{\cal M}$  is acyclic than it is a simplex with vertices  $Y_1, \ldots.Y_n$  and  $X_{[n]'}$. We claim that  $\tilde{Y}:=\{Y_1, \ldots Y_n\}$  is a facet of  $\tilde{\cal M}$.

To prove the claim we  use the axioms for maximal vectors (topes) of an oriented matroid \cite{IS2}  
to conclude that the signed set  $\tilde{X}[0]:=( \{X_{A'}\}_{A\subseteq [n]},\emptyset) $,  complementary of  $\tilde{Y}$, is a positive cocircuit     of  $\tilde{\cal M}$.

Let  $\tilde{\cal T}$  denote the family of maximal covectors (topes) of the adjoint  $\tilde{\cal M}$  of  $\bf O_n$, consider the family  $\tilde{\cal P}$  of principal maximal covectors defined by:  $\tilde{\cal P}:=\{ \tilde{X}[e_1]\circ ...\circ \tilde{X}[e_n]:\  e_i\in \{i,i'\},\ \forall i\in[n]\} $. Taking the restriction of (all the vectors of)  $\tilde{\cal P}$  to  $\tilde{Y}$  we conclude:
$$\tilde{\cal P}(\tilde{Y})= \{+,-\}^{\tilde{Y}}\leqno{(1)}$$ 
This equality implies  (\cite {IS2} )  that  $\tilde{Y}$  is an independent set of  $\tilde{\cal M}$.
  
On the other hand, the restriction to  ${\cal D} _1\setminus \tilde{Y}$  of  every vector of $\tilde{\cal P}$   is  the same vector:
$$\tilde{X}[0]:=( \{X_{A'}\}_{A\subseteq [n]},\emptyset) \leqno{(2)}$$
From  $(1)$  and  $(2)$  we conclude that  $\tilde{X}[0]\circ \tilde{\cal T}\subseteq \tilde{\cal T}$  implying that  $\tilde{X}[0]$  is a positive covector of  $\tilde{\cal M}$  with  $\tilde{Y}$  as the complementary  flat.  Since $\tilde{Y}$  is independent   $r(\tilde{Y})=r(\tilde{\cal M})-1$  implying   $\tilde {Y}$  is actually a facet of  $\tilde{\cal M}$. 

{\parindent=0cm
{\bf 3)} We concluded in 2)  that the positive sign vector  $\tilde{X}[0]$  is a cocircuit of   $\tilde{\cal M}$ and therefore, by Theorem A,  a localization in  $\bf O_n$, with respect to  ${\cal D}_1$.  
}
The single element extension  ${\bf O_n}\cup 0$  of  ${\bf O_n}$  obtained by localizing the new element  $0$  according to  $\tilde{X}[0]$  contains as signed cocircuits the signed-vectors of the following set:
$$\overline{\cal D}_1:=\{ X_{A'}\cup 0=( A'\cup ([n]\setminus A)\cup 0, \emptyset)\}_{A\subseteq [n]} \cup \{Y_i=(i,i')\}_{i\in [n]}.$$ 
We leave to the reader the verification that  for every hyperpline  $L$  of  $\bf O_n$  the hyperplane of  ${\bf O_n}\cup 0$  spanned by  $L\cup 0$  always contains a hyperplane of  $\bf O_n$, implying that the family   $\overline{\cal D}$ is actually the family of signed cocircuits  of the extension   ${\bf O_n}\cup 0$, proving 3).\vskip 2mm

{\parindent=0cm
{\bf 4)} For every  $i\in [n]$,  the isomorphisms   ${\bf O_n}/i \simeq {\bf O_n}/i' \simeq  {\bf O_{n-1}}\cup 0$  follow directly from the definition of  the cross-polytope  $\bf O_n$, the details are left to the reader.
}

The fact that $\tilde{\cal M}(\tilde{H}_i)$ is an adjoint of   ${\bf O_n}/i$  with respect to  ${\cal D}_1/i$  is a consequence of the hypothesis that  $\tilde{\cal M}$  is an adjoint of  ${\bf O_n}$  with respect to  ${\cal D}_1$  and the fact that the family of signed cocircuits of   ${\bf O_n}/i$  is  ${\cal D}/i:=\{ X=(X^+,X^-)\in {\cal D} : \ i\in X^0\}= {\cal D}_1/i\cup -{\cal D}_1/i= \tilde{H}_i\cup -\tilde{H}_i$.

The case  $\tilde{\cal M}(\tilde{H}_{i'})$  is similiar.

\vskip 1cm

\begin{Theo}

The canonical oriented adjoint   $\tilde{\bf O}_{\bf n}$  of the cross-polytope  $\bf O_n$  described in Proposition 2.1   is, up to isomorphism, the unique \underline{realizable} oriented adjoint, with respect to  ${\cal D}_1$, of  ${\bf O_n}$. 

\end{Theo}

{\parindent=0cm

{\bf Proof.}  The proof is by induction on  $n$. For  $n=2$  the result is trivially true. Assume that the result is true for  $k<n$  and consider a realizable  adjoint   $\tilde{\cal M}=\tilde{\cal M}({\cal D}_1)$   of   $\bf O_n$.   By definition  ${\cal D}_1= \tilde{H}_n\cup \tilde{H}_{n'}\cup  Y_n$  and by induction hypothesis, and Proposition 2.2-3) and 4),  the restriction of  $\cal M$  to each one of the hyperplanes  $\tilde{H}_n=\{ X_{A'\cup n'}\}_{A\subseteq [n-1]}\cup (\tilde{Y}\setminus Y_n)$  and  $\tilde{H}_{n'}=\{ X_{A'}\}_{A\subseteq [n-1]} \cup (\tilde{Y}\setminus Y_n)$  is isomorphic to the canonical adjoint of  $\bf O_{n-1}$.
}

By Proposition 2.2.1) and 2)  $\cal M$  is acyclic with vertices  $Y_1, \ldots Y_n$  and  $X_{[n]'}$. Let  $f: {\cal D}_1\longrightarrow \R^n$  be an affine representation of  $\cal M$. Consider  $f(Y_i)=y_i$,  $f(X_{[n]'})=b$  and  $ f(X_\emptyset)=c$.    $f(\tilde{H}_n)=B(Y_1, \ldots, Y_{n-1}, X_{[n]'})$  is contained in the affine hyperplane  $H= aff(y_1\ldots, y_{n-1},b)$  of  $\R ^n$. The point   $f(X_n)\in relint conv(y_1,\ldots,y_{n-1},b)$   and   $B=(y_1, \ldots,y_n,b)$  is an affine basis of  $\R^n$. 

From proposition  2.2.1 we know that  $c$  must be a point in the relative interior of the line segment  $y_n f(X_n)$, therefore a point in the interior of  $conv(B)$  and   $f(\tilde{H}_{n'}$  is contained in the affine  hyperplane  $H'=aff(y_1, \ldots, y_n,c)$.  Then, for each  $i\in [n]$  the circuit  $\tilde{C}(i, \emptyset)$  defines the point   $f(X_i)$,  as the intersection point  $f(X_i)= aff(y_i, c)\cap conv  ( \{y_k\}_{k\in [n]\setminus i}\cup b) =  aff(y_i, c)\cap conv( f(\tilde{H}_{i} ))$.  By Proposition 2.1  $ f(\tilde{H}_{i})$ is completely determined by  the affine base  $B\setminus y_i$  and  by  $c_i$, yelding that for every  $A\subset [n]$    $f_(X_{A'})=F_{[n]\setminus A}\cap S_A$  and   $\tilde{\cal M}=B(Y_1, \ldots, Y_{n}, X_{[n]'})= \tilde{\bf O}_{\bf n}$.
\vskip 5mm

\begin{section}{Cubes, adjoints and extension spaces of the cross-polytope}

We start by recalling  the definitions of cubes and of oriented cubes from  \cite{IS1}:\vskip 5mm

{\parindent=0cm
{\bf Notation.} We consider   $C^n:=\{-1,1\}^{[n]}$. 

For every  $A\subseteq [n]$,  $v_A$ denotes the vector of  $C^n$  defined by:  $v_A(i)=-1\ for \ i\in A$  and  $v_A(i)=1\ for \ i\notin A$. It is sometimes useful to write  $v_A=(-1_A, 1_{[n]\setminus A})$.

Given  $B\subseteq [n]$  and  $v\in C^n$,  $_{-B}v$  denotes the vector of  $C^n$  whose coordinates are the opposite of those in $v$  in the entries indexed by  $B$, i. e.  $_{-B}v(i)=-v(i)\ for\ i\in B$  and   $_{-B}v(i)=v(i)\ for\ i\notin B$.

Notice that  $_{-B}v_A=  v_{A\Delta B}$  where  $A\Delta B:= A\setminus B\cup B\setminus A$. 

Given a subset  $V\subseteq  C^n$  we denote by  $_{-A}V$  the subset of  $C^n$  defined by:  $_{-A}V:=\{_{-A}v\in C^n:\ v\in V\}$.}
\vskip 5mm

The {\it real cube} is the oriented matroid   ${\bf Q_n}=Aff(C^n)$  of the affine dependencies of  $C^n$  over the reals.

Notice that the circuits with smallest number of elements of the real cube are the {\it signed rectangles}  $R(v; A,B)=(\{ v,_{-(A\uplus B)}v\}, \{_{-A}v,_{-B}v\}),\ with  \ A\cap B=\emptyset$. We denote by  $\cal R$  the family of all the signed rectangles of  $\bf Q_n$.  

A {\it facet} of the real cube is a hyperplane  $H_{i+}:=\{v\in C^n: \ v(i)=1\}$   or   $H_{i-}:=\{v\in C^n: \ v(i)=-1\}$   and a {\it skew-facet} is a hyperplane  $H_{ij+}:=\{v\in C^n: \ v(j)=v(i)\}$  or  $H_{ij-}:=\{v\in C^n: \ v(j)=-v(i)\}$. We denote by  $\cal F$  the signed cocircuits of the real cube complementary to the facets   and by  $\overline{\cal F}$  the signed cocircuits complementary of the facets and to the skew-facets.
  $\overline{\cal F}$  is the collection of  signed cocircuits of  $\bf Q_n$  with smallest cardinality.
  \vskip 1cm

{\parindent=0cm
\begin{Def} { \rm( Cubes and  Oriented Cubes \cite{IS1})} 

A {\it cube over $C^n$ }  is a matroid  $M=M(C^n)$  satisfying the following two conditions:

C1) every (unsigned) rectangle of  $C^n$  is a circuit of  $M$

C2) every facet and skew facet of  $C^n$  is an (unsigned) cocircuit and a hyperplane of  $C^n$. 

Cubes may or may not be orientable. Orientable cubes have a canonical orientation  ${\cal M}(C^n)$  that satisfies the following equivalent conditions:

C1) all the signed rectangles of  $\bf Q_n$  are circuits of  $\cal M$.

C2) all the signed cocircuits of  $\bf Q_n$  that are facets and skew-facets are signed cocircuits of  $\cal M$.

An {\it oriented cube} is a cube  isomorphic to an oriented cube  over  $C^n$, ${\cal M}(C^n)$,  canonically oriented.  
\end{Def}
\vskip 2mm

{\bf Remark.} Note that to conclude that an oriented matroid  ${\cal M}(C^n)$  is an oriented cube one only has to verify the following conditions:

C1) every signed rectangle is a circuit. 

C'2) the facets are hyperplanes of   $\cal M$. 

\vskip 5mm

\begin{Prop}
Let  $\tilde{\cal M}=\tilde{\cal M}({\cal D}_1)$  be an adjoint of the cross-polytope   ${\bf O_n}$  with respect to  ${\cal D}_1:=\{Y_i \}_{i\in [n]} \cup \{ X_{A'} \}_{ A\subseteq [n]}$. Let  ${\cal D}_1^+$  be the subset of  ${\cal D}_1$  defined by  ${\cal D}_1^+:=\{ X_{A'} \}_{ A\subseteq [n]}$  of  ${\cal D}_1$.

The restriction   $\tilde{\cal M}({\cal D}_1^+)$  of  $\tilde{\cal M}$  to  ${\cal D}_1^+$   is an oriented cube.
\end{Prop}

{\bf Proof.} Let  $\tilde{\cal M}':=\tilde{\cal M}({\cal D}_1^+)$  be the restriction of an adjoint of  ${\bf O_n}$ to  ${\cal D}_1^+$. 
We identify  ${\cal D}_1^+$  with  $C^n$  via the bijection that to every  $X_{A'},\ A\subseteq[n]$  assigns the vector  $v_A\in C^n$ and prove that ${\cal M}({\cal D}_1^+)$  satisfies the  conditions  $C1)$  and  $C'2)$  of  the above Remark.

C'2) is a direct consequence of the definitions of adjoint and of restriction. Notice that for every  $i\in [n]$  the restriction  of the signed cocircuits  $\tilde{X}[i]$  and  $\tilde{X}[i']$  of  $\tilde{\cal M}$  to  ${\cal D}_1^+$  are the positive cocircuits  of $\tilde{\cal M}'$:

$\tilde{X}[i]({\cal D}_1^+)=(\{X_{A'}\}_{ A\subseteq[n]\setminus i}, \emptyset)$  and $\tilde{X}[i]({\cal D}_1^+)=(\{X_{A'\cup i'}\}_{A\subseteq[n]\setminus i}, \emptyset)$ 

whose complementary sets in  ${\cal D}_1^+$, are the  LV-facets of  $\tilde{\cal M}'$ :

$H_{i^+}:= \tilde{H_i}({\cal D}_1^+)=\tilde{X}[i]^0({\cal D}_1^+)$  and   $H_{i-}:=\tilde{H}_{i'}({\cal D}_1^+)=\tilde{X}[i]^0({\cal D}_1^+)$.
}

The lattice of faces of  $\tilde{\cal M}'$  is therefore isomorphic to the lattice of faces of  $\bf Q_n$.

{\parindent=0cm
C1)  We prove that every signed rectangle   $(\{{X},_{-(A\uplus B)}{X}\}, \{_{-A}{X}, _{-B}{X}\})$  is a circuit of   $\tilde{\cal M}'$,  where for every   $X\in {\cal D}_1^+$, $X=X_{C'}$  for some  $C\subset[n]$,  and every  $D\subset[n] $,  $_{-D}X=X_{(C\Delta D)'}$. The proof is by induction on  $n$.
}

For $n=2$ it is trivial,  in this case  $\tilde{\cal M}$  must be the oriented matroid represented in Figure 1  and clearly  $\tilde{\cal M}'$  is a rectangle, a cube of rank $3$.

\begin{figure}[h] 
\begin{center}
\includegraphics[scale=0.70]{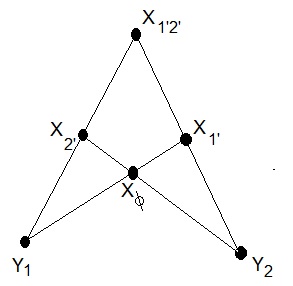}
\end{center}
\vskip -5mm
\caption{ $\tilde{O}_2$, the unique adjoint of  $\bf O_2$}
\end{figure}

Assume that C1)  is true for every adjoint  $\tilde{\cal M}$  of  $\bf O_{n-1}$.

Consider an adjoint  $\tilde{\cal M}$  of  $\bf O_{n}$, the restrictions  $\tilde{\cal M}(\tilde{H}_n)$  and  $\tilde{\cal M}(\tilde{H}_{n'})$  are  isomorphic to adjoints of  $\bf O_{n-1}$ (Proposition 2.2.4)), therefore all the rectangles contained in the facets  $\tilde{H}_{n+}$  and of  $\tilde{H}_{n-}$  of  $\tilde{\cal M}'$  are rectangles of  $\tilde{\cal M}'$  and we are left with proving that  the signed rectangles of the form   $R(X; A,B\uplus n)=(\{X, _{-(A\uplus B\uplus n)}X\}, \{_{-A}X, _{-(B\cup n)}X\})$, with  $X\in {\cal D}_1^+$, $A\uplus B\subseteq [n-1]$  are signed circuits of  $\tilde{\cal M}'$. This will be done by induction on the cardinality of  $B$  assuming, without loss of generality, that  $n\in X^+$, i.e. $X$  belongs to the facet  $\tilde{H}_{n+}$  of  $\tilde{\cal M}'$
 
When  $|B|=0$, {\it every  rectangle  $R(X; A, n)=(\{X, _{(A\cup n}X\}, \{_{-A}X, _{- n}X\})$  with  $A\subset[n-1]$  is a circuit of  $\tilde{\cal M}'$. }

In fact, by Proposition 2.2.1), since  $n\in X^+$  the signed circuits  $C=(\{Y_n, _{-n}X\}, X)$  and  $D=(\{Y_n, _{-(A\cup n)}X\}, _{-A}X)$  are circuits of  $\tilde{\cal M}$. Eliminating  $Y_n$  between  $-C$  and  $D$  we conclude that  $R(X; A, n)$  is a circuit of  $\tilde{\cal M}'$. 

Assume that  $R(X; A, B\cup n)$  is a circuit of  $\tilde{\cal M}'$  for  $|B|=k$  and consider  a rectangle  $R_B:= R(X; A, B\cup n)$  with  $|B|=k+1$.

Consider  $B_b=B\setminus b$, for some  $b\in B$. The rectangle  $R(X; A, B_b\cup n)=(\{X, _{-(A\cup B_b\cup n)}X\},\{ _{-A}X,_{-(B_b\cup n)}X\})$  is a circuit of  $\tilde{\cal M}$. 

\hskip -1cm Consider the rectangle   $S= (\{_{-(B_b\cup n)}X, _{-(A\cup B\cup n)}X\}, \{_{-(B\cup n)}X\}, _{-(A\cup B_b\cup n)}X\})$  contained in the hyperplane  $H_{n'}$  of  $\tilde{\cal M}$  and therefore a rectangle  of   $\tilde{\cal M}'$. Eliminating  $_{-(B_b\cup n)}X$  between  $R(X; A, B_b\cup n)$  and  $S$  we conclude that $R_B$  must be a signed rectangle of  $\tilde{\cal M}'$, proving the Lemma.

\vskip 0.5cm
The next Theorem describes completely the relation between oriented cubes and adjoints of the cross-polytope.
\vskip 5mm

\begin{Theo}
Consider   $\tilde{\cal O}_n$  the class of all oriented adjoints of the cross-polytope  $\bf O_n$  with respect to  ${\cal D}_1:=\{Y_i \}_{i\in [n]} \cup \{ X_{A'} \}_{ A\subseteq [n]}$  and  $\tilde{\cal C}_n$  the class of all oriented cubes over  ${\cal D}_1^+:=\{ X_{A'} \}_{ A\subseteq [n]}$.

The correspondence  $f: \tilde{\cal O}_n\longrightarrow \tilde{\cal C}_n$  defined by  $f(\tilde{\cal M})=\tilde{\cal M}({\cal D}_1^+)$,  that assigns to every oriented adjoint  $\tilde{\cal M}\in \tilde{\cal O}_n$  its restriction  to  ${\cal D}_1^+$,  is a bijection.
\end{Theo}\vskip 5mm

The proof of Theorem 3.1 requires proving that every oriented cube can  be extended, in a unique way, to an oriented matroid isomorphic to an adjoint of the cross-polytope. In order to do so we need some preliminary lemmas.
\vskip 5mm

{\parindent=0cm
\begin{Lema} {\rm (A property of flats and covectors of cubes)}

Consider an oriented cube  ${\cal M}={\cal M}(C^n)$ and a flat  $F$  of  $\cal M$.

For every  $i\in [n]$  let  $_{-i}F:=\{_{-i}v: \ v\in F\}$. Then,

1) Either  $F\cap\ _{-i}F=\emptyset$  or   $F=\ _{-i}F$. 

2) Let  $V=(V^+,V^-)$  be a vector of  $\cal M$  such that  $V^0=F$. 

\hskip 5mm (i) if  $F\cap\ _{-i}F=\emptyset$,  define  $F+:= F\cap H_{i^+}$  and  $F-:= F\cap H_{i^-}$, the following holds:  either   $_{-i}(F+)\subseteq V^+$  and  $_{-i}(F-)\subseteq V^-$  or   $_{-i}(F+)\subseteq V^-$  and  $_{-i}(F-)\subseteq V^+$.

\hskip 5mm (ii) if  $F=_{-i}F$  then  $V(H_{i^-})= _{-i}V(H_{i^+})$.
\end{Lema}

{\bf Proof.} 1) We prove that if  $F\cap _{-i}F\not=\emptyset$  then  $F=_{-i}F$: assume that there is some element  $v\in F\cap\  _{-i}F$, i.e.  $v,_{-i}v\in F$.  If for some  $A\subseteq [n]\setminus i$,  $_{-A}v\in F$,  also  $_{-(A\cup i)}v$  must belong to  $F$  since the rectangle   $(\{v,_{-(A\cup i)}v\},\{_{- i}v, _{-A}v\})$  is a circuit of  $\cal M$. 1) is proved.

2) Assume that  $V=(V^+,V^-)$  is a covector of  $\cal M$  complementary to the flat  $F$. In both cases we use orthogonality between covectors and the rectangles of the cube  $\cal M$  to obtain the result:

(i) if  $F\cap\ _{-i}F=\emptyset$, then for every  $v,v'\in F+$  orthogonality  of  $V$  with the rectangle  $(\{v,_{-i}v'\}, \{_{-i}v,v'\})$  implies that  $_{-i}v,_{-i}v'\in   _{-i}(F+)$  must have the same sign in  $V$. 
}
For every  $v\in F+$  and  $w\in F-,  \ w=_{-(A\cup i)}v$  for some  $\ A\subseteq [n]\setminus i$,  orthogonality  of  $V$  with the rectangle  $(\{v,w\}, \{_{-i}v,_{-i}w\})$  implies that  $_{-i}v\in _{-i}(F+)$  and  $_{-i}w\in _{-i}(F-)$  must have different signs in  $V$. The result then, follows.

\vskip 5mm
{\parindent=0cm
\begin{Def}{\rm Flats parallel and anti-parallel to  $i$}

{\it Let  $F$  be a flat of an oriented cube  ${\cal M}={\cal M}(C^n)$.

The flat  $F$  is {parallel to $i$} if it satisfies the condition:  $F=_{-i}F$.

The flat  $F$  is {anti-parallel to $i$} if it satisfies the condition:  $F\cap _{-i}F=\emptyset$.}
\end{Def}
}
\vskip 5mm 

\begin{Lema} {\rm (Extensions by the point at infinity of a class of parallel edges of a cube)}  Consider an oriented cube  ${\cal M}=(C^n,{\cal D} )$,  ${\cal D}$  its family of signed cocircuits.
For every  $i\in[n]$  define 

${\cal Y}_i:=\{ Y=(Y^+,Y^-)\in {\cal D}: Y^0$ {\it is a hyperplane parallel to  $i$} $\}$

${\cal X}_i:=\{ X=(X^+,X^-)\in {\cal D}: X^0=(F+)\uplus (F-)$ {\it is a hyperplane anti-parallel to  $i$ and} $\  _{-i}(F-)\subseteq X^+ \ and  _{-i}(F+)\subseteq X^- \}$, where  $(F+):=X^0\cap H_{i^+}$  and  $(F-):=X^0\cap H_{i^-}$

Then  

1. $({\cal X}_i,{\cal Y}_i)$  is a localization in $\cal M$  defining a single element extension  ${\cal M}\cup y_i$  of  $\cal M$.

2. $\cal M$  and its  single element extension  ${\cal M} \cup y_i$  have the same number of hyperplanes and signed cocircuits.

The new element  $y_i$ will be called the point at infinity of direction  $i$.
 \end{Lema}

{\parindent=0cm
{\bf Proof.}  From the definitions of  ${\cal X}_i$  and  ${\cal Y}_i$  it is clear that   ${\cal D}= {\cal X}_i\uplus -{\cal X}_i\uplus {\cal Y}_i$.

For 1.  we have to prove that for every hyperpline  $L$  of  ${\cal M}$,  $({\cal X}_i\uplus {\cal Y}_i)/L$  is a localization in the rank 2 realizable oriented matroid  ${\cal M}/L$.  

Let  $L$  be a hyperline of  ${\cal M}$.

{\it Case 1)  If  $L$  is parallel to $i$} then,  by Lemma 3.1, every hyperplane containing  $L$  is also parallel to  $i$  and  in this case   $({\cal X}_i\uplus {\cal Y}_i)/L=  {\cal Y}_i/L=  {\cal D}/L $, and  $y_i$  is a loop in  ${\cal M}/L$.

{\it Case 2) If $L$  is anti-parallel to $i$} then  $L\cap\ _{-i}L=\emptyset$. Consider  $v\in L$  and  the hyperplane $H= cl_{\cal M}(L\cup _{-i}v),\ v\in L$. Since for every  $w\in L\setminus v$, the set  $\{v,w,_{-i}v,_{-i}w\}$  is a circuit of  $\cal M$  we conclude that   $H=L\cup\ _{-i}L$ and therefore  $H$  is the unique hyperplane parallel to  $i$  that contains  $L$. Moreover, for every hyperplane anti-parallel to  $i$, containing  $L$,  its complementary signed cocircuit   $X=(X^+,X^-)$  in  ${\cal X}_i$  is the signed cocircuit that satisfies the condition  $_{-i}(L-)\subseteq X^+$  and  $_{-i}(L+)\subseteq X^-$, where   $(L+):= L\cap H_{i^+}$  and  $(L-):= L\cap H_{i^-}$. The localization of the new point  $y_i$  in  ${\cal M}/L$  is represented in the next figure, clearly a localization in  ${\cal M}/L$.
 }
 
\begin{figure}[h] 
\begin{center}
\includegraphics[scale=0.70]{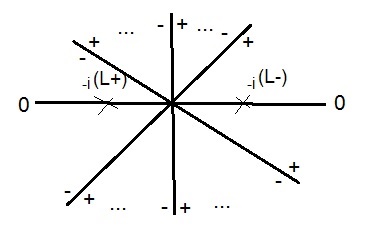}
\end{center}
\vskip -5mm
\caption{ Localization of $y_i$  in ${\cal M}/L$}
\end{figure}

Now consider  2. By definition of the extension  ${\cal M}\cup y_i$,  for every  $X=(X^+,X^-)\in {\cal X}_i$  the  signed set  $X+{y_i}=(X^+\cup y_i,X^-)$   is a signed cocircuit of  ${\cal M}\cup y_i$  and its  complement of  $(X+{y_i})^0=X^0$ is a hyperplane of both  $\cal M$  and   ${\cal M} \cup y_i$.  For every pair of opposite signed cocircuits  $\pm Y=(Y^+,Y^-)\in {\cal Y}_i$, its complement, $Y^0$, in $C^n$  is  a hyperplane of  $\cal M$  and in  $C^n\cup y_i$, $Y^0\cup y_i$,  is a hyperplane of  ${\cal M} \cup y_i$. 

We claim that  ${\cal M}\cup y_i$ has no other hyperplanes.  Any other hyperplane of  ${\cal M}\cup y_i$  would be of the form  $L\cup y_i$  for some hyperline  $L$  of  ${\cal M}$. There is no such hyperplane if  $L$  is parallel to  $i$, because in this case  $y_i\in L$. On the other hand if  $L$  is  anti-parallel to $i$  the unique hyperplane of  ${\cal M}\cup y_i$  that contains  $y_i$ is the hyperplane  $H=L\cup _{-i}L \cup y_i$  and  $H\setminus y_i$  is an hyperplane of  $\cal M$. The result is proved. \vskip 5mm

{\parindent=0cm
{\bf Remark.}  Notice that by definition of the localization  $({\cal X}_i,{\cal Y}_i)$  the extension  ${\cal M}\cup y_i$  contains as hyperplanes all the sets   $H_{j^+}\cup y_i$,  $H_{j^-}\cup y_i$  for  $j\not=i$  implying that for every  $v\in C^n$  the set  $\{v,_{-i}v, y_i\}$  is a circuit and a flat of  ${\cal M}\cup y_i$. Moreover,  supposing that  $v\in H_{i^+}$  the signature of this circuit must be   $\pm C$  where  $C=( \{_{-i}v, y_i\},\{v\})$  since for every hyperline  $L$  intersection of a facet  $H_{j^{\epsilon}}$  with  $H_i^+$ ,  $y_i$ is placed in  every hyperplane containing  $L$, opposite to  $_{-i}L= H_{j^\epsilon} \cap H_{i^-}$.
\vskip 1cm 

{\bf Proof of Theorem 3.1.} In Proposition 3.2 we proved that the restriction to   ${\cal D}_1^+$  of any adjoint of  ${\bf O_n}$,  with respect to  ${\cal D}_1$  is an oriented cube. 
}

In order to prove the theorem we consider an oriented cube,  ${\cal M}={\cal M}(C^n)$ and prove that it can be extended, in a unique way, to an oriented matroid isomorphic to an adjoint of the cross-polytope.

By Lemma 3.2  ${\cal M}$,  admits a single element extension by a new point $y_i$, the point at infinity in direction  $i$.

To conclude the proof of the Theorem we have to prove that the  $n$  single-element  extensions  $y_1, \ldots, y_n$  of the cube $\cal M $  are compatible, resulting in a unique oriented matroid  ${\cal M}(C^n\cup \{y_1, \ldots, y_n\})$  that is isomorphic to an adjoint of  ${\bf O_n}$.

Since no single-element extension by a point  $y_i$  creates new hyperplanes (Lemma 3.2), we can perform successively, in any order, the  $n$  elementary extensions by elements   $y_i$. Note, that the description of the localizations defining the extension  of the oriented matroid   ${\cal M}(C^n\cup \{y_{i_1},y_{i_2}, \ldots , y_{i_k}\}),\ k<n$  by  $y_{i_{k+1}}$  is exactly the same localization that describes the single element extension of  ${\cal M}(C^n)$   by  $y_{i_{k+1}}$.

When the  $n-th$  extension by an element $y_{i_n}$  is completed, independenly of the order, we obtain the same matroid  ${\cal M}(C^n\cup\{y_1, \ldots, y_n\})$. 

Consider the bijection   $f:C^n\cup\{y_1,\ldots,y_n\} \longrightarrow {\cal D}_1$  defined  for every $i\in [n]$  by  $f(y_i)= Y_i=(\{i\},\{i'\})$, and for every  $A\subseteq [n]$  by  $f(v_A)=X_{A'}$. By construction, the image  $f({\cal M}(C^n \cup\{y_1, \ldots, y_n\}))$  of the extended matroid   ${\cal M}(C^n\cup\{y_1, \ldots, y_n\})$  is an adjoint of the cross-polytope  ${\bf O_n}$  with respect to  ${\cal D}_1$.
\vfill\eject

\begin{Theo}
If  ${\cal M}(C^n)$  is a \underline{realizable} oriented cube then  ${\cal M}(C^n)\simeq {\bf Q_n}$ 
\end{Theo}
\vskip 5mm

To prove the theorem we prove that every realization   ${\cal A}ff(W)$  of  $\cal M$  can be extended to a realizable adjoint of the cross-polytope  $\bf O_n$. Then Theorem 2.1 ensures that  ${\cal M}\simeq {\bf Q_n}$.

We use the next two lemmas (actually from lemma 3.4 we only need (i)).
\vskip 5mm

\begin{Lema}{\rm ( in Theorem 1 of \cite{IS1} )}
Let   ${\cal M}={\cal M}(C^n)$ be an oriented cube. Then for every  $i\in [n]$  the reorientation of  ${\cal M}$  obtained by reversing signs on a facet  $H_{i^-}$  is an oriented cube  $_{-H_{i^-}}{\cal M}=_{-H_{i^-}}{\cal M}(C^n)$  whose facets are the facets  $H_{i-}$, $H_{i+}$   and the  skew-facets  $H_{ij+}, H_{ij-}$, $j\in [n]\setminus i$, of  $\cal M$ .
\end{Lema}

{\parindent=0cm
{\bf Proof.} The proof is straightforward from the definitions and is left to the reader.
\vskip 5mm

\begin{Lema}
Let  ${\cal M}={\cal M}(C^n)$  be a realizable oriented cube of rank  $n+1$.  Let  $W\subset \R^n$  be a realization of  $\cal M$ such that  ${\cal M}\simeq {\mathcal A}ff(W)$  via the bijection  $\sigma:  C^n\longrightarrow W$,  $\sigma (v_A)=w_A$  for every  $A\subseteq [n]$. Then the following two conditions are satisfied for  $n\geq 2$:

(i) All the line segments  $w_Aw_{[n]\setminus A}$  with  $A\subseteq [n]$  meet in the same point  $O\in \R^n$. The point  $O$ is a point in the interior of the polytope  $conv(W)$  and  will  be called the \underline{center of (the realization)  $W$}.

(ii) For  $n\geq 3$, let  $O_{i^+}$,  $O_{i^-}\in \R$  denote the centers of the facets   $H_{i^+}, H_{i^-}$  of   $W$,  let  ${S}:=\{ O_{i^+}$,  $O_{i^-}:\ i\in [n]\}$.  Then the oriented matroid  $\mathcal Aff(S)$  is a cross-polytope - a polar of the the polytope  $conv(W)$.

\end{Lema}
\vskip 2mm

{\bf Notation.} To simplify notation we use the same notation  $H_{i^+}$, $H_{i-}$,  $H_{ij+}$,  $H_{ij-}$  for the facets and skew-facets of  ${\cal M}$  and of its realization.  
By an {\it opposite pair of points of the realization  $W$}  is meant any pair of points the form  $w_A$, $w_{[n]\setminus A}$.\vskip 2mm

{\bf Proof.}  (i) Let  $w_A$, $w_B$  and  $w_C$, $A,B,C\subseteq[n-1]$,  be three distinct points of the facet  $H_{n^+}$  of  ${\cal M}_W:={\cal A}ff(W)$. The opposite points  $w_{[n]\setminus A}$, $w_{[n]\setminus B}$  and  $w_{[n]\setminus C}$  are points of the facet $H_{n^-}$. The signed rectangle  $C(A,B)=(\{w_A,w_{[n]\setminus A}\}, \{w_B,w_{[n]\setminus B}\})$  is a circuit of  ${\mathcal M}_W$  and therefore  the four points are the  vertices of a quadrilateral  whose diagonals  $w_Aw_{[n]\setminus A}$  and  $w_Bw_{[n]\setminus B}$  intersect in a point   $O\in relint(w_Aw_{[n]\setminus A})\cap relint (w_Bw_{[n]\setminus B})$.  
}

Similarly, considering the circuits  $C(A,C)= (\{w_A,w_{[n]\setminus A}\}, \{w_C,w_{[n]\setminus C}\})$  and  $C(B,C)= (\{w_B,w_{[n]\setminus B}\}, \{w_C,w_{[n]\setminus C}\})$
we conclude that the diagonal  $ w_Cw_{[n]\setminus C}$  of both quadrilaterals must intersect the other diagonal. Let  $O_A$  be the intersection point of the diagonals of  $C(A,C)$  and  $O_B$  the intersection point of the diagonals of  $C(B,C)$ (Figure 3)

\begin{figure}[h] 
\begin{center}
\includegraphics[scale=0.70]{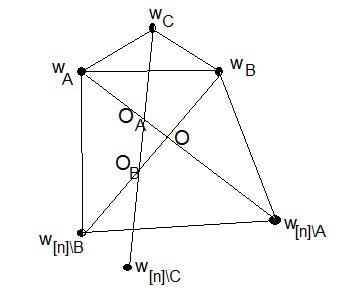}
\end{center}
\vskip -5mm
\caption{ }
\end{figure}

To conclude the proof of   (i)  we claim that :  $O_A=O_B=O.$

In fact, if  $O_A\not=O_B$  then the line segment  $ w_Cw_{[n]\setminus C}$  would be contained in the affine plane  $\pi:= aff (C(A,B))$. On the other hand, since no cube contains three colinear points , the plane  $\pi$  would be the affine plane  spanned by the points  $w_A,w_B,w_C$. But, since  $n\geq 3$, this affine plane  is contained  in the facet  $H_{n^+}$  and therefore it cannot contain points of  $H_{n^-}$, otherwise  $rank({\cal M})<n+1$. But this contradicts the assumption that  $C(A,B)\subset \pi$. We conclude then that  $O_A=O_B$  and, consequently, it must be  $O_A=O_B=O$.

By definition $O$ is in the relative interior of every line segment connecting opposite vertices of  $W$. It is clearly a point of the convex polytope  $conv(W)$  which is in none of its facets, therefore  $O\in int (conv(W))$.

{\parindent=0cm
(ii) We prove that  $O\in relint (O_{n^+}O_{n^-})$, where  $O_{n^+}$, $O_{n^-}$  are the centers of the facets  $H_{n+}$, $H_{n-}$  of  ${\cal M}_W$. The other cases are similar.
}

Let  $w_A, w_B$  be two distinct points of the facet  $H_{n^+}$  such that  $B\not= A, [n-1]\setminus A$. This is possible  because  $n\geq 3$. 

Consider the circuits   $C=(\{w_Aw_{[n-1]\setminus A}\}, \{w_Bw_{[n-1]\setminus B}\})$  contained in  $H_{n^+}$  and  $C'= (\{w_Aw_{[n]\setminus A}\}, \{w_Bw_{[n]\setminus B}\})$  contained in  $H_{n^-}$.  By definition  $O_{n^+}:= relint (w_{A}w_{[n-1]\setminus A})\cap relint(w_{B}w_{[n-1]\setminus B}) $  and 
$O_{n^-}:=relint (w_{A\cup n}w_{[n]\setminus A})\cap relint (w_{B\cup n}w_{[n]\setminus B}) $. 

Consider now the signed rectangles of  $W$  defined as:    

\hskip -1cm $D=(\{w_A,w_{[n]\setminus A}\}, \{w_{[n-1]\setminus A},w_{A\cup n}\})$  and   $D' =(\{w_B,w_{[n]\setminus B}\}, \{w_{[n-1]\setminus B},w_{B\cup n}\})$. 

Denote by  $conv(D), conv(D')$  the quadrilaterals of  $\R^n$  whose vertices are the vertices of  the circuits  $D$  and  $D'$.
The line segment  $O_{n^+}O_{n^-}= conv(D)\cap conv(D')$ (figure 4).

\begin{figure}[h] 
\begin{center}
\includegraphics[scale=0.70]{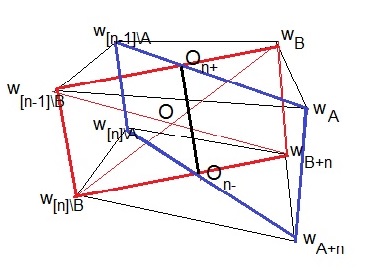}
\end{center}
\vskip -5mm
\caption{ }
\end{figure}

On the other hand since, by definition, the center   $O$  is the intersection of the diagonals of the  $D$  and  $D'$, it follows that  $O\in relint O_{n^+}O_{n^-}$.

In a similar way we conclude that   $O\in relint O_{i^+}O_{i^-}$, for every  $i\in[n]$,  implying that  $\forall   \ i<j\leq n$   $C_{ij}=(\{O_{i^+},O_{i^-}\}, \{O_{j^+},O_{j^-}\})$  is a circuit of  the matroid  $\mathcal Aff(S)$, $S=\{O_{i^+},O_{i^-}: i\in [n]\}$. Since the points  $O_{i^+}$  and  $O_{i^-}$  lie in the relative interior of opposite facets of the  convex polytope  $conv(W)$  that has dimension  $n$, the  $n$  vectors   $O_{i^-}O_{i^+}$  form a basis of  $\R^n$implying that for every   $i\in [n]$,  $S\setminus \{O_{i^+},O_{i^-}\}$  must be a hyperplane of  $\mathcal Aff(S)$. This oriented matroid is an n-cross-polytope. 
\vskip 1cm

{\parindent=0cm
{\bf Proof of Theorem 3.2} }

Let  ${\cal M}_W={\cal A}ff(W)$  be an affine realization in  $\R^n$  of a (realizable) oriented cube  ${\cal M}={\cal M}(C^n)$  via the correspondence  $v_A\longrightarrow w_A\in \R^n, \ \forall A\subseteq [n]$.  By Lemma 3.3 the reorientation  of  $\cal M$,  $_{-H_{i-}}{\cal M}$, obtained by reversing signs on the facet  $H_{i-}$  is an oriented cube whose facets are the  2n-hyperplanes  of  ${\cal M}$:  $H_{i^+}$,  $H_{i^-}$,  $H_{ij+}$,  $H_{ij-}$, for every  $j\in [n]\setminus i$.

Considering  $\R^n$  as a subset of the real projective space  $\mathbb P_n(\R)$, representations of the oriented cube   $_{-H_{i-}}{\cal M}$  can be obtained as images  $f_i(W)$  of  $W$  under suitable projective transformation  $f_i:\mathbb P_n(\R)\longrightarrow \mathbb P_n(\R)$  admissible for  $W$. In every such representation  $f_i(W)$  of  $_{-H_{i-}}{\cal M}$, the images, $f_i(w_A)$, $f_i(w_{A\cup i})$  of the extremities of the edges  $w_Aw_{A\cup i}$  of  $W$   are opposite points of   $f_i(W)$. By lemma 3.4  all of the line segments  $f_i(w_A)f_i(w_{A\cup i})$  meet in the center  $C_i$  of  $f_i(W)$. Therefore  all the edges of  ${\cal A}ff(W)$   meet in the  point  $f_i^{-1}(C_i)\in \mathbb P_n(\R)$.  This shows that every realization  ${\cal A}ff(W)$  of  $\cal M$  can be extended to a realizable adjoint of the cross-polytope. By Theorem 2.1 there is a unique realizable adjoint of the cross-polytope. By Theorem 3.1  ${\cal A}ff(W)$  must be isomorphic to  the real affine cube  ${\cal A}ff(C^n)$.\vskip 2mm

{\parindent=0cm
{\bf Remark.}  Actually a suitable projective transformation   $f_i$  can be constructed in the following way, using the notation of Lemma 3.3:
identify, as usual, $\R^n$  with  $\R^n\times {1}\subseteq \R^{n+1}$  in such a way that the center  $O_W$  of  $W$  has coordinates  $({\bf 0},1)$. Let $\pi_i$  be an affine plane of  $\R^{n+1}$   cutting the line segment  $O_WO_{i^+}$ in such a way that it separates the points of the facet  $H_{i^-}$  and the origin of  $\R^n$   from  the points of  $H_{i^+}$. The projection  $f_i(W)$  from the origin  $O\in \R^{n+1}$  of  $W$ in the plane  $\pi$  is a representation of  $_{-H_{i-}}{\cal M}$.
Let  $C_i$  be the center of  $f_i(W)$. The point   $f_i^{-1}(C_i)$  of the real projective space  $\mathbb P_n(\R)$ is the meeting point of all the edges  $w_Aw_{A\cup i}$  of  $W$. 
\vskip 5mm

The next Theorem 3.3. is a nonoriented version of Theorem 3.2.
\vskip 5mm
 
\begin{Theo}
Let  ${\cal M}(C^n)$  be an oriented cube whose underlying matroid  $\underline{\cal M}(C^n)$  is  {realizable}. Then  $\underline{\cal M}(C^n)\simeq \underline {\bf Q_n}$. 
\end{Theo}
\vskip 5mm

The proof of this Theorem is an easy consequence of Theorem 3.2 and of the following Theorem B  from \cite{IS1}:
\vskip 5mm

{\bf Theorem B}( {\rm Corollary 1.3)  of \cite{IS1}})
{\it Let  $M(C^n)$  be an oriented cube whose underlying matroid is realizable. Then every class of orientations of  $\cal M$ contains a unique acyclic orientation whose signed rectangles and signed cocircuits complementary of the facets and skew-facets  are signed as in  ${\cal A}ff(C^n)$.} 
\vskip 5mm

{\bf Proof of theorem 3.3.}  Let  $W\subseteq \R^n$  be a realization of the underlying matroid. By lemma 3.4 we may assume without loss of generality that  ${\cal A}ff(W)$  is a realizable oriented cube. Theorem 3.2  
then ensures that  $ {\cal A}ff(W)\simeq {\bf Q_n}$ and the result follows.
}
\end{section}

\begin{section}{Final Remarks}

From Theorem 3.2 we know that  if an oriented cube other than the real cube exists, it must be non-realizable. Moreover, from Theorem 3.3, we know that it may be non-realizable either because it is a non-realizable orientation of the real affine cube, or because its underlying matroid is a cube that is non-representable over  $\R$.

The following conjectures play a fundamental role in this context: \vskip 2mm

{\parindent=0cm
{\bf Conjecture 1.} (M. Las Vergnas \cite{LV3}): {\it The real affine cube  $\bf Q_n$  has a unique class of orientations.} Equivalently, {\it the canonical adjoint   ${\tilde{\bf O}}_{\bf n}$   of the cross-polytope  $\bf O_n$  has a unique class of orientations}.\vskip 2mm

{\bf Conjecture 2.} (I. P. Silva \cite{IS1}): {\it The real affine cube  $\bf Q_n$  is the unique oriented cube.} Equivalently, {\it the canonical adjoint    ${\tilde{\bf O}}_{\bf n}$  of the cross-polytope is the unique adjoint of the cross-polytope  $\bf O_n$.}\vskip 2mm

}

Several  constructions and examples of (unoriented) cubes non-representable over  $\R$ are given    \cite{IS1}.  All the known examples are non-orientable matroids.  New families of minor-minimal excluded minors for orientability based on extensions of the cross-polytope  are described in  \cite{IS3}. \vskip 2mm

\end{section}

\end{document}